\documentclass[11pt,draftclsnofoot,onecolumn]{IEEEtran}   

\usepackage[pdftex]{graphicx}
\usepackage[cmex10]{amsmath}
\usepackage{cancel}
\usepackage{amssymb}
\usepackage{color}
\usepackage{epstopdf}
\usepackage{hyperref}
\usepackage{epsfig}
\usepackage{bm}
\usepackage{subfig}
\usepackage{ifthen}
\usepackage{bibunits}
\usepackage{longtable}
\usepackage{ulem}
\usepackage{algorithm}
\usepackage{algorithmic}

\DeclareSymbolFont{AMSb}{U}{msb}{m}{n}
\DeclareMathSymbol{\R}{\mathalpha}{AMSb}{"52}
\DeclareMathSymbol{\LZ}{\mathalpha}{AMSb}{"5A}
\DeclareMathSymbol{\N}{\mathbin}{AMSb}{"4E}

\newtheorem{theorem}{Theorem}[section]

\newboolean{FIGSINBACK}
\setboolean{FIGSINBACK}{FALSE}   

\begin{document}

\title{Reducing Basis Mismatch in Harmonic Signal Recovery via Alternating Convex Search}

\author{Jonathan M. Nichols*, Albert K. Oh, Rebecca M. Willett \\ 

\thanks{J. M. Nichols* (jonathan.nichols@nrl.navy.mil) is with the U. S. Naval Research Laboratory, Washington, D.C. 20375. A. K. Oh (akoh2@wisc.edu) and R. M. Willett (rmwillett@wisc.edu) are with the Department of Electrical and Computer Engineering, University of Wisconsin-Madison, Madison, WI 53706. Code associated with this paper is available at: \url{https://sites.google.com/site/albertktoh/software/acs}.}
}

\maketitle

\begin{abstract}

The theory behind compressive sampling pre-supposes that a given sequence of observations may be exactly 
represented by a linear combination of a small number of basis vectors. In practice, however, even small deviations from an exact signal model can result in dramatic increases in estimation error; this is the so-called ``basis mismatch'' problem.  This work provides one possible solution to this 
problem in the form of an iterative, biconvex search algorithm.  The approach uses standard $\ell_1$-minimization to find the 
signal model coefficients followed by a maximum likelihood estimate of the signal model. The algorithm is illustrated on harmonic signals of varying sparsity and outperforms the current state-of-the-art.

\end{abstract}


\section{Introduction \label{sec:intro} }

One of the most fundamental tasks in signal processing is the estimation of the coefficients associated with a given signal model.  This problem frequently takes the form of a linear inverse problem ${\bf y=\Psi x+\bm{\eta} }$ where ${\bf y}\in\R^N$ are the observed data, ${\bf x}\in\R^N$ are the model coefficients and $\bm{\eta} \in\R^N$ is a ``noise'' vector, chosen from some joint probability distribution, usually Gaussian.  A good model adheres to the principle of parsimony, allowing $S<N$ non-zero coefficients to accurately describe the observations, in which case we say the model is $S$-sparse \cite{Donoho:09}.  For real-valued signals, a common signal model ${\bf \Psi}\in\R^{N\times N}$ consists of sine and cosine vectors (i.e. the Fourier basis), in which case the Discrete Fourier Transform (DFT) of the data yields the estimates ${\bf \hat{x}}$. 

The problem becomes more complex, however, if we seek to estimate a high-fidelity model using low-fidelity data.  
The compressed sensing (CS) framework introduced by Donoho \cite{Donoho:06}, Candes {\it et al.} \cite{Candes:06a}, and others considers the sampling model
\begin{align}
{\bf y}={\bf A\Psi x}+{\bm\eta}
\label{eqn:1}
\end{align}
where the matrix ${\bf A}\in \R^{M\times N}$ projects the high-dimensional signal ${\bf z=\Psi x}$ onto the (potentially low-dimensional) observations ${\bf y}\in\R^M$, and $\bm{\eta} \in\R^M$.  For the undersampled case, where $M<N$,  conditions for accurate estimation are well-established and mandate that $S<\rho M$ where the fraction 
$\rho\in (0,1)$ depends on the degree of undersampling, $\delta=M/N$ \cite{Donoho:10}.  The greater the degree of undersampling, the smaller $\rho$ must be for accurate solutions.  If $\rho$ is small enough, and ${\bf A}$ is chosen appropriately (e.g., entries are random draws from some probability distribution), then (\ref{eqn:1}) can be solved via
\begin{align}
\hat{\bf x}= \underset{\bf x}{\arg \min} \ \|{\bf y-A\Psi x}\|_2^2+\lambda\|{\bf x}\|_1,
\label{eqn:L1}
\end{align}
where $\lambda$ is a positive constant that penalizes non-sparse solutions and where $\|{\bf x}\|_1:=\sum_{k=1}^{N} |x_k|$ is large.   The initial theoretical work on this estimator \cite{Candes:05, Donoho:06} has given rise to new approaches to super-resolution \cite{Candes:12}, analog-to-digital conversion \cite{Tropp:10, Nichols:11, Yenduri:12}, and imaging \cite{Romberg:08,lustig2007sparse,willett2014sparsity,harmany2011spatio}.

The key to accurate estimation is the sparsifying transformation matrix ${\bf \Psi}$, chosen so that only $S<M$ elements of the solution ${\bf x}$ are non-zero.  For example, the sum of $S$ harmonics can be described by exactly $S$ complex Fourier coefficients ($2S$ real coefficients) of the aforementioned Fourier basis.  However, because we observe discrete signals, only certain 
frequencies can be represented exactly.  The Fourier basis for an $N$-point signal (sampled at unit intervals) is comprised of sine and cosine vectors at discrete frequencies $f_k=(k-1)/N,~k=1,\ldots,  N/2 + 1$.  So long as the signal frequencies match {\it exactly} those in this set, 
the representation (Fourier basis) indeed yields a $S$-sparse representation.  However, if the signal frequencies lie off this frequency grid (e.g., $f=(k+0.5)/N$ for $k$ integer), some number $J>S$ frequencies will be required in the representation. Unfortunately, $J$ may be quite a bit larger than $S$, precluding accurate reconstruction of ${\bf x}$ with only $M = S/\rho$ measurements \cite{Nichols:13}.

This is the crux of the so-called ``basis mismatch'' problem, whereby even a good signal model can yield a poor reconstruction due to seemingly small differences between basis vectors assumed by the reconstruction algorithm and a similar set yielding a far sparser representation of the signal.  This problem has received recent theoretical
\cite{Calderbank:11,Gilbert:12,Tang:12} and experimental \cite{Ekanadham:11,Nichols:13,Rao:12} treatment and solutions are still under development.  One straightforward approach to the basis mismatch problem is to simply oversample the frequency space (e.g., let ${\bf \Psi} \in \R^{N \times QN}$ contain sinusoids with frequencies 1/QN apart instead of 1/N apart for a large integer Q).  Though a higher resolution frequency discretization can shrink errors, it comes at the cost of an 
increasingly underdetermined problem with stronger coherence between the vectors in the signal model, both of which work against the conditions needed for successful recovery.

Another approach is to treat both the model coefficients {\it and} the associated frequencies as unknowns to be solved. Boufounos {\it et al.} developed an approach that isolates the unknown frequencies in separate, non-overlapping bins and then solves for their 
location and amplitude, though the method is restricted to specific sampling strategies \cite{Gilbert:12}.  For more general sensing matrices, Tang {\it et al.} have shown that both frequency locations and amplitudes can be estimated by 
solving the constrained ``atomic norm minimization'' problem via semi-definite programming (SDP) \cite{Tang:12}.  Rao {\it et al.} also formulate the problem as one of atomic norm minimization but propose a greedy ``forward-backward'' (GFB) algorithm as the solution \cite{Rao:12}.  The forward-backward approach lacks the theoretical guarantees of the SDP, however it runs in faster time for larger problems.

In this work, we present an Alternating Convex Search (ACS) algorithm \cite{Gorski:07} that significantly reduces the influence of basis mismatch on the recovery of a linear combination of sinusoids.  The approach iterates between a constrained minimization using the familiar $\ell_1$-norm, and a component-wise minimization 
involving the model vectors themselves.  Rather than relying on a random gridding of the frequency space, as was done in \cite{Rao:12}, we estimate the frequency vectors that minimize the mean squared error between the data and model at each step. The method is computationally fast relative to SDP, although slightly less accurate, and appears more inline with the approach described in \cite{Rao:12} as well as the work by Valley and Shaw \cite{Valley:13}.

\section{Algorithm Description \label{sec:algorithm}}

Consider the real harmonic continuous time signal $z(\tau)=\sum_{k=1}^{S} x_k \exp(i2\pi f_k \tau)$ with frequencies $f_k\in [0,1/2]$ and complex coefficients given by $x_k$.  We can accurately model $N$ values of this signal as ${\bf z=\Psi x}$ taking ${\bf \Psi}$ to be the $N\times N$ discrete Fourier basis.  As our sampling model, we consider the aformentioned compressed sensing (CS) framework and write ${\bf y}={\bf Az}+{\bm\eta}={\bf A\Psi x}+{\bm\eta}$ where the entries of the $M\times N$ matrix ${\bf A}$ are chosen in accordance with CS theory to be independent draws from some probability distribution \cite{Donoho:10}.  As we have mentioned, the estimator (\ref{eqn:L1}) will degrade (possibly fail) if the signal frequencies are not chosen from the set $f_k=(k-1)/N,~k=1,\ldots, N/2+1$ (see e.g., \cite{Nichols:13}). 

Our solution is to treat the frequencies in the signal model as unknowns to be estimated.  Rather than restrict ourselves to the Fourier basis, we consider the more general overcomplete harmonic ``dictionary'' where each of the columns are referred to as ``atoms'' \cite{Donoho:01}.  We let $Q$ be a factor specifying the degree of overcompleteness, chosen so that $QN/2$ is an integer, and form the $N\times QN$ real-valued matrix ${\bf \Psi}_{\bm\theta}$ where the entries are given by

\small
\begin{align}
&\left({\bf\Psi}_{\bm\theta}\right)_{n,k}=\nonumber \\
&~\left\{\begin{array}{r@{,~}l}
		\sqrt{\frac{2}{N}}\cos\left(2\pi n\left(\frac{k-1}{QN}+\theta_k\right)\right) & k=1,\ldots,  \frac{QN}{2} \\
		-\sqrt{\frac{2}{N}}\sin\left(2\pi n\left(\frac{N-k}{QN}-\theta_{N-k+1}\right) \right) & k=\frac{QN}{2}+1,\ldots,  QN
		\end{array}\right.\nonumber \\
~&n=0,\ldots,  N-1
\label{eqn:basisparam}
\end{align}
\normalsize
where each of the frequency parameters $\theta_k,~k=1,\ldots,  QN/2$ can take on values in the range $[-\frac{1}{2QN},\frac{1}{2QN}]$. Although such a representation provides continuous coverage over the frequency space, for $S\ge 6$ (three or more harmonics) the $f_k$ must be spaced at least $1/(QN)$ apart in order to be resolved.  Note also that though there are $QN$ indices, there are only $QN/2$ elements of $\bm{\theta}$ corresponding to the $QN/2$ unique frequencies in the parameterized basis set. This has implications for the implementation of our approach, as any change to the ``positive'' frequencies must be accompanied by an equal and opposite change to the ``negative'' frequencies (sidestepping complex exponentials this way allows us to leverage off-the-shelf sparse optimization tools for real-valued signals and bases.).
The end result is $QN/2$ unknown frequency perturbations that must be estimated.

Clearly if $Q=1$ and all of the 
${\bm\theta}$ are zero, then (\ref{eqn:basisparam}) defines the standard Fourier basis.  However, in the general case we modify (\ref{eqn:L1}) and 
propose solving the following minimization problem
\begin{align}
\left(\hat{\bf x},{\hat{\bm\theta}}\right)&=\underset{{\bf x},{\bm\theta}}{\arg \min} \ f({\bf x},{\bm\theta}) \equiv \|{\bf y-A\Psi}_{\bm\theta}{\bf x}\|_2^2+\lambda\|{\bf x}\|_1,
\label{eqn:cost}
\end{align}
which simply states that we would like to find both a sparse signal model ${\bf \Psi}_{\bm\theta}$ by estimating the ${\bm\theta}_k$ and the associated model coefficients ${\bf x}$.  

In the first step we solve the familiar $\ell_2 \mbox{-} \ell_1$ problem
\begin{align}
{\bf \hat{x}}^{(t+1)}=\underset{\bf x}{\arg \min} \ \|{\bf y-A\Psi}_{{\hat{\bm\theta}}^{(t)}}{\bf x}\|^2_2+\lambda\|{\bf x}\|_1,
\label{eqn:cost1}
\end{align}
holding perturbations ${\bm \theta}$ fixed. This is a convex optimization problem that is known to yield unique sparse solutions provided that the problem size is chosen in accordance with Figure 5 of  \cite{Donoho:10}; the solution will be sparse with only $S<QN$ non-zero coefficients.

The next step is to hold the model coefficients ${\bf x}$ fixed, and find a better signal model by solving for the frequency perturbations of non-zero coefficients. Note that we declare a coefficient to be ``non-zero'' if its magnitude exceeds a threshold $\kappa$ and we capture the associated indices in the set $\mathcal{S}^{(t+1)} \equiv \{\ell:|x_\ell^{(t+1)}| \ge \kappa\}$. In principle, one could set $\kappa=0$, in which case the algorithm must minimize over all $QN/2-1$ unknown perturbations ${\bm\theta}$; the cost is a significant increase in CPU time.  On the other hand, too large a threshold risks ``missing'' a frequency component that requires adjusting. The updated signal model is found via
\begin{align}
{\hat{\bm\theta}}^{(t+1)}_{\mathcal{S}^{(t)}} &=\underset{{\bm \theta}_{\mathcal{S}^{(t)}} }{\arg \min} \ \|{\bf 
y-A\Psi}_{\bm \theta}{\bf \hat{x}}^{(t+1)}\|_2^2.
\label{eqn:cost2}
\end{align}
As implied by the notation, rather than solving for all unknown frequency parameters, we solve only for ${\bm\theta}_{\mathcal{S}^{(t)}}$, which are the indices of ${\bm \theta}$ in the set $\mathcal{S}^{(t)}$. To solve for \eqref{eqn:cost2}, we leverage that each of the atoms in our collection are nearly linearly independent, so rather than searching the $S$-dimensional frequency space, we divide the problem into a series of $S$ one-dimensional searches; this simplification was also noted in \cite{Stoica:89}.  
Thus, for each non-zero frequency index $\ell \in \mathcal{S}^{(t)}$ we solve the one-dimensional problem
\begin{align}
\hat{\theta}_\ell^{(t+1)} =  \underset{\theta_\ell}{\arg \min} \ \|{\bf y-A\Psi}_{{\bm \theta}}{\bf \hat{x}}^{(t+1)}\|_2^2
\label{eqn:cost3}
\end{align}
while holding each of the other parameters fixed to their current values.

We repeat these two stages of the algorithm until a convergence criteria is met, as described in Algorithm \ref{alg:ACS_alg}.

\begin{algorithm}
\caption{ACS Algorithm}
\label{alg:ACS_alg}
\begin{algorithmic} 

\STATE Input signal $\bf y$ and sensing matrix $\bf A$. Select $Q$.
\STATE Set $\alpha = 0.1$, $\beta = 0.1$, and $TOL = 10^{-5}$
\STATE Set $t \leftarrow 0 $
\STATE Initialize $\hat{{\bm \theta}}^{(t)} = 0$

\REPEAT
\STATE Compute ${\bf \Psi}_{\hat{\bm \theta}^{(t)}}$ via \eqref{eqn:basisparam}
\STATE $\lambda = \alpha \| ({\bf A} {\bf \Psi}_{{\hat{\bm\theta}}^{(t)}})^T {\bf y} \|_{\infty}$
\STATE $\hat{{\bf x}}^{(t+1)} \leftarrow$ solution of \eqref{eqn:cost1}
\STATE $\kappa = \beta \| {\hat{\bf x}}^{(t+1)} \|_2$
\STATE Create $\mathcal{S}^{(t+1)} \equiv \{\ell:|\hat{x}_\ell^{(t+1)}| \ge \kappa\}$

\FOR {$\ell \in \mathcal{S}^{(t+1)}$}
\STATE ${\hat{\theta}_\ell^{(t+1)}} \leftarrow$ solution of \eqref{eqn:cost3}
\ENDFOR

\STATE $\hat{f}^{(t+1)} \triangleq f(\hat{{\bf x}}^{(t+1)}, \hat{{\bm \theta}}^{(t+1)})$
\STATE $t \leftarrow t+1$
\UNTIL  {$| (\hat{f}^{(t+1)} - \hat{f}^{(t)})/\hat{f}^{(t)} | < TOL$}

\STATE Output $\hat{{\bf x}} = \hat{{\bf x}}^{(t)}$, $\hat{{\bm \theta}} = \hat{{\bm \theta}}^{(t)}$

\end{algorithmic}
\end{algorithm}

In the supplementary material, we show that the function (\ref{eqn:cost3}) is convex on the frequency range $\theta_k\in [-1/(2QN),1/(2QN)]$ so long as $Q\ge 3/2$ and the true frequencies are separated by at least $1/QN$.  Hence, our algorithm is appropriately described as an Alternating Convex Search (ACS) \cite{Gorski:07}. A benefit to this strategy is that we may make use of Theorem 4.5 in \cite{Gorski:07}, which states that for a biconvex function $f({\bf x},{\bm\theta})$, the sequence $f(\{ {\bf x}^{(i)},{\bm\theta^{(i)}} \})_{i\in \N}$ generated by ACS converges monotonically.  Although convergence can only be guaranteed for an overcomplete dictionary, we found that setting $Q=1$ is sufficient and even advantageous in practice (when the true tones are sufficiently well-separated) as the resulting problem is not as underdetermined.  Nonetheless, results will be shown for both $Q=1$ and $Q=3/2$.

An interesting feature of this algorithm is that it will sometimes shift two frequency atoms toward the same value.  For example, if the signal frequency lies near $k/2QN$ for some k, the algorithm may set $\theta_{k-1} = 1/2QN$ and $\theta_k = -1/2QN$, in which case the associated coefficients, $\hat{x}_{k-1}$ and $\hat{x}_k$, will sum to the true value (i.e., the signal energy is split between two atoms lying at or near the same frequency line).  From a sparsity standpoint this is somewhat undesirable as a single coefficient is preferred.  However, the resulting model is perfectly valid and produces accurate results at the cost of only two additional real coefficients.  This effect will be clearly seen in the results that follow (see Figure \ref{fig:coeff3}, second frequency line).

\section{Numerical Experiments \label{sec:experiment}}

To study the proposed algorithm, we construct a harmonic signal with $S/2$ non-zero complex 
($S$ real) coefficients.  Specifically we set, for $i=1,\ldots, N$,
\begin{align}
z(\tau_i)=\sqrt{2/N}\sum_{s=1}^{S/2} \cos\left(2\pi \left(f_s+\theta_s\right) \tau_i+\varphi_s\right),
\label{eqn:sigmodel}
\end{align}
where each frequency is given a uniformly distributed phase $\varphi_s\in [0,2\pi)$. The frequencies themselves are chosen $f_k=k/QN+\theta_k$ where $\theta_k\sim Uniform(-1/(2QN),1/(2QN))$ and \\ $k\sim Uniform(\{0,1,\ldots,QN/2-1\})$ (without replacement). The sampling interval was assumed to be $\tau_{i+1}-\tau_i=1$.

We consider $M=128$ projections as the measured signal ${\bf y}$.  From these projections, we seek to infer the $N=256$ values of $\bm{x}$ that define the underlying signal ${\bf z}={\bf \Psi x}$. Our method easily scales to much larger values of $N$ and $M$, but we consider smaller problem sizes here to facilitate comparisons with the much slower SDP.  Each entry of the ${\bf A}$ matrix is assumed known and chosen independently from a standard normal distribution.  The sequence of noise values $\bm{\eta}$ are also chosen independently from a normal distribution with standard deviation $\sigma=\left(\frac{1}{N}\sum_{i=1}^{N} z(\tau_i)^2\right)^{1/2}/(10^{(\mathrm{SNR \ dB}/20)})$, where $\mathrm{SNR \ dB}$ is the signal-to-noise ratio (log base 10 scale); we chose $\mathrm{SNR \ dB}=40$.  

\begin{figure}[t]
\begin{center}
\includegraphics[width=3in]{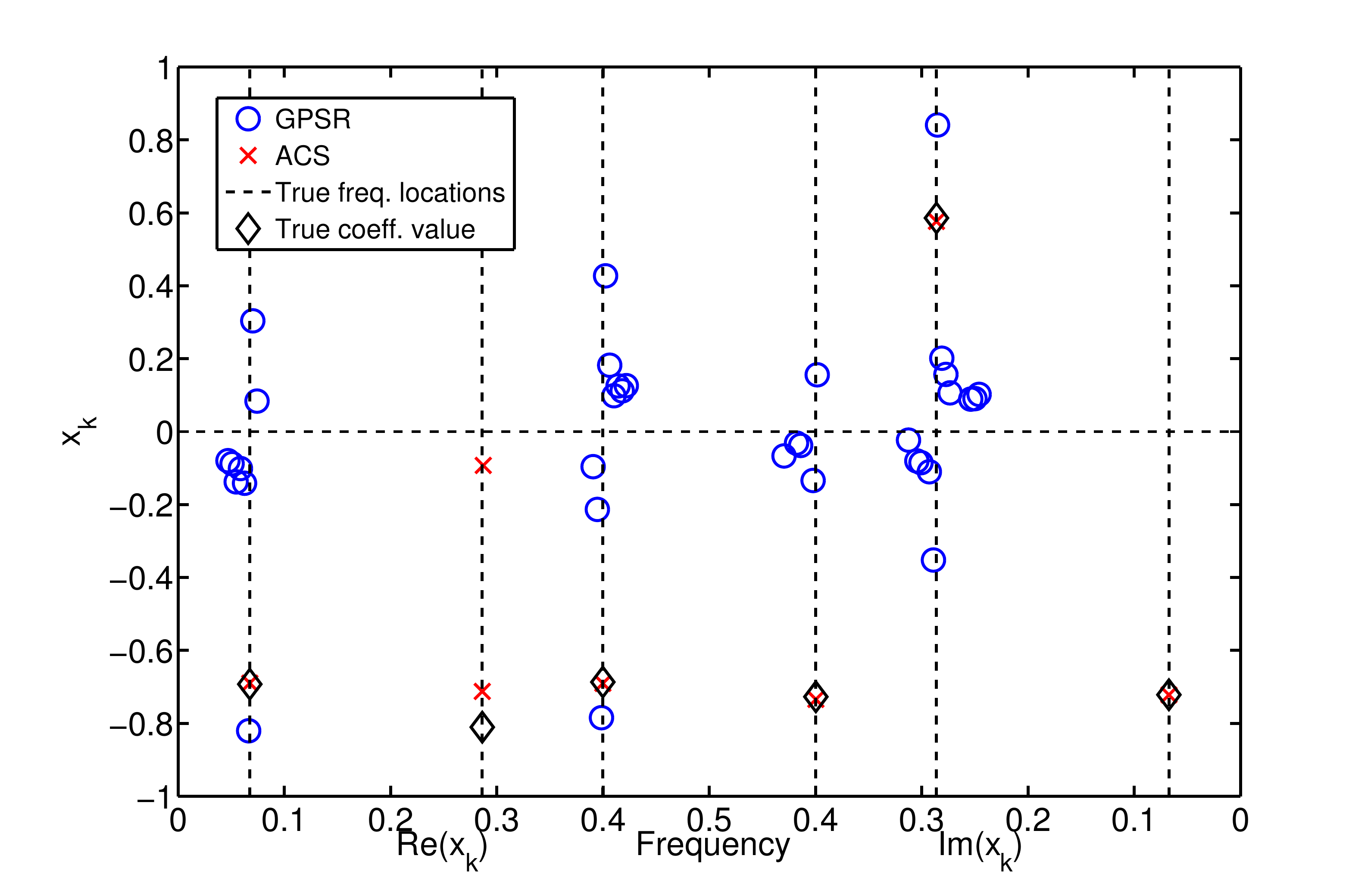}
 \caption{Real and imaginary components of the identified coefficients $x_k$ given in the frequency ranges $[0,0.5]$ and $(0.5,0)$, respectively. Shown are true coefficient values and those identified by standard GPSR and the ACS algorithm.} 
   \label{fig:coeff3}
\end{center}
\vspace{-3ex}
\end{figure}

We note that the Gradient Projections for Sparse Reconstruction (GPSR) algorithm \cite{gpsr} was used to solve the $\ell_2 \mbox{-} \ell_1$ optimization problem in \eqref{eqn:cost1}. It was also used to solve \eqref{eqn:L1} using the standard Fourier basis as a baseline comparison. Note that the GPSR algorithm for both setups used identical termination criteria and we set the recommended $\lambda = \alpha \| ({\bf A} {\bf \Psi})^T {\bf y} \|_{\infty}$ value. As a first example we set $S=6$; Figure \ref{fig:coeff3} shows the coefficients obtained from ACS, the standard GPSR implementation, and true coefficient values. ACS clearly produces a sparser model and correctly identifies the true frequency locations and corresponding magnitudes.  

\newlength{\figwidth}

\setlength{\figwidth}{0.4\linewidth}

\subsection{Comparison to Overcomplete Fourier Dictionary}

In order to address this basis mismatch problem, perhaps the simplest solution is to oversample the frequencies by creating an overcomplete Fourier dictionary $\bm{\Psi}_{OC} \in \R^{N \times QN}$, which is not parameterized with respect to a frequency perturbation, but rather includes a finer discretization of the frequency components. The atoms in this dictionary are identical to those defined in (\ref{eqn:basisparam}), excluding the frequency perturbations. Our discretized signal ${\bf z}$ is then modeled to be sparse in this dictionary by writing ${\bf z} = {\bf\Psi}_{OC} {\bf w}$, where the coefficients ${\bf w} \in \R^{QN}$ can be solved via (\ref{eqn:L1}).  This approach is also 
implemented with GPSR and will be referred to as OC-GPSR.

Figure \ref{fig:oc_comp} compares the reconstruction performance of ACS and OC-GPSR for two measures of model error and algorithm execution time.  Specifically, we analyze the normalized RMSE, defined as $\|{\bf z}-\hat{{\bf z}}\|_2^2/\|{\bf z}\|_2^2$, as well as a measure of the identified model support.  Denote as $\mathcal{J}_i\equiv \{j: ~|\hat{f}_j-f_i|<\epsilon\}$ the set of indices for which the identified tones are within $\epsilon$ Hz of the $i^{th}$ tone and compute
\begin{align}
\mathrm{err}=\sum_{i=1}^S |x_i-\sum_{j\in\mathcal{J}_i} \hat{x}_j|.
\label{eqn:err1}
\end{align} 
This measure captures the ability of the approach to correctly identify the support of the true model (see e.g., \cite{Recht:14}).  In this work we always use $\epsilon=1/(5QN)$, i.e., 1/5 of a frequency bin, as the threshold for closeness.

Perhaps the biggest advantage to the ACS approach is correct identification of the frequency support (small $\mathrm{err}$).  The OC-GPSR will typically place many frequencies near the true frequency, resulting in a low RMSE.  However the resulting model is not a parsimonious description of the data and is clearly over-fitting the model.  In terms of execution time the two methods are comparable for $Q>8$, however for mildly overcomplete dictionaries the OC-GPSR approach is a good deal faster which may be essential for some applications.

\setlength{\figwidth}{0.4\linewidth}
\begin{figure}[t]
\begin{center}
\subfloat{\includegraphics[width=\figwidth]{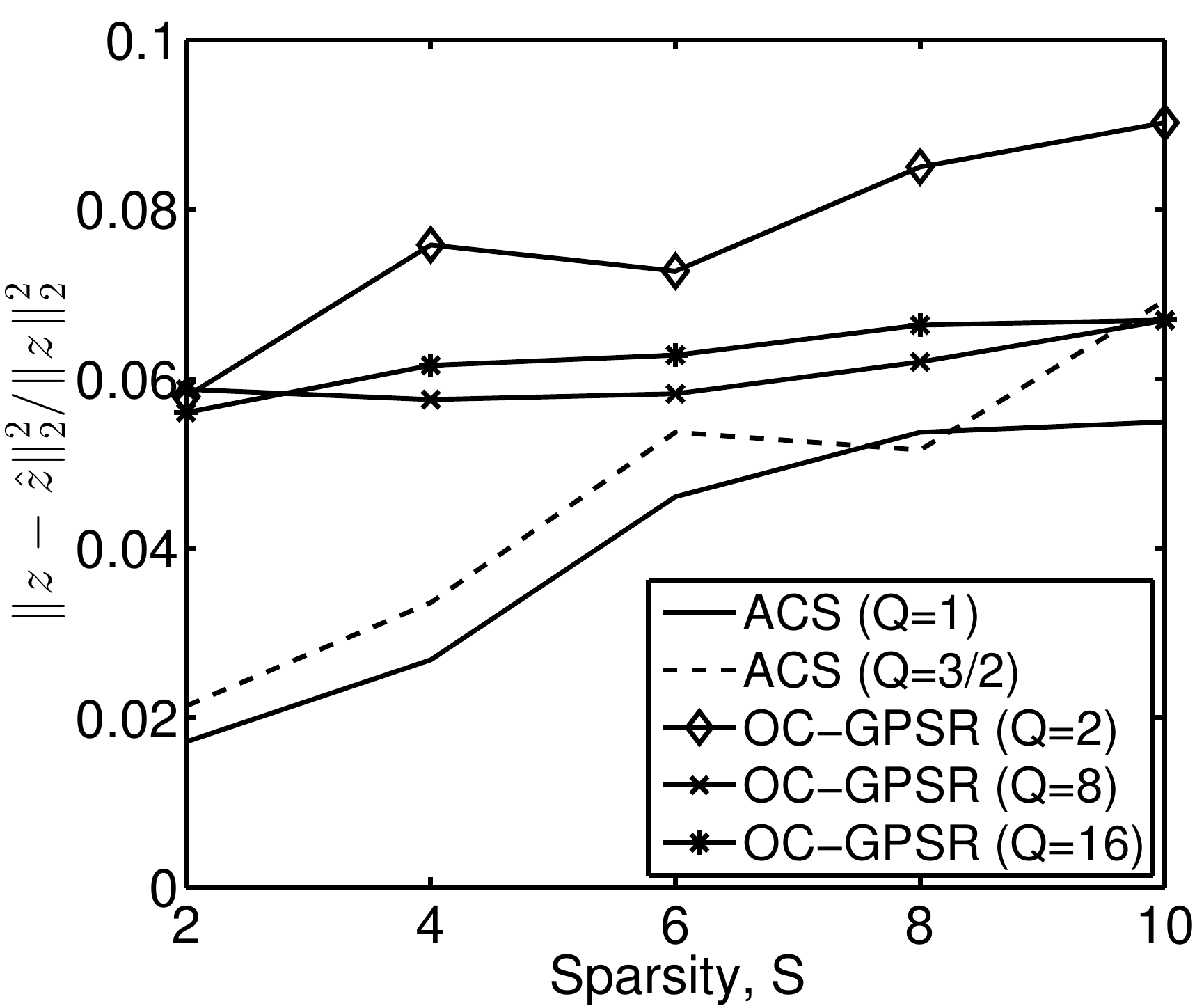}} ~
\subfloat{\includegraphics[width=\figwidth]{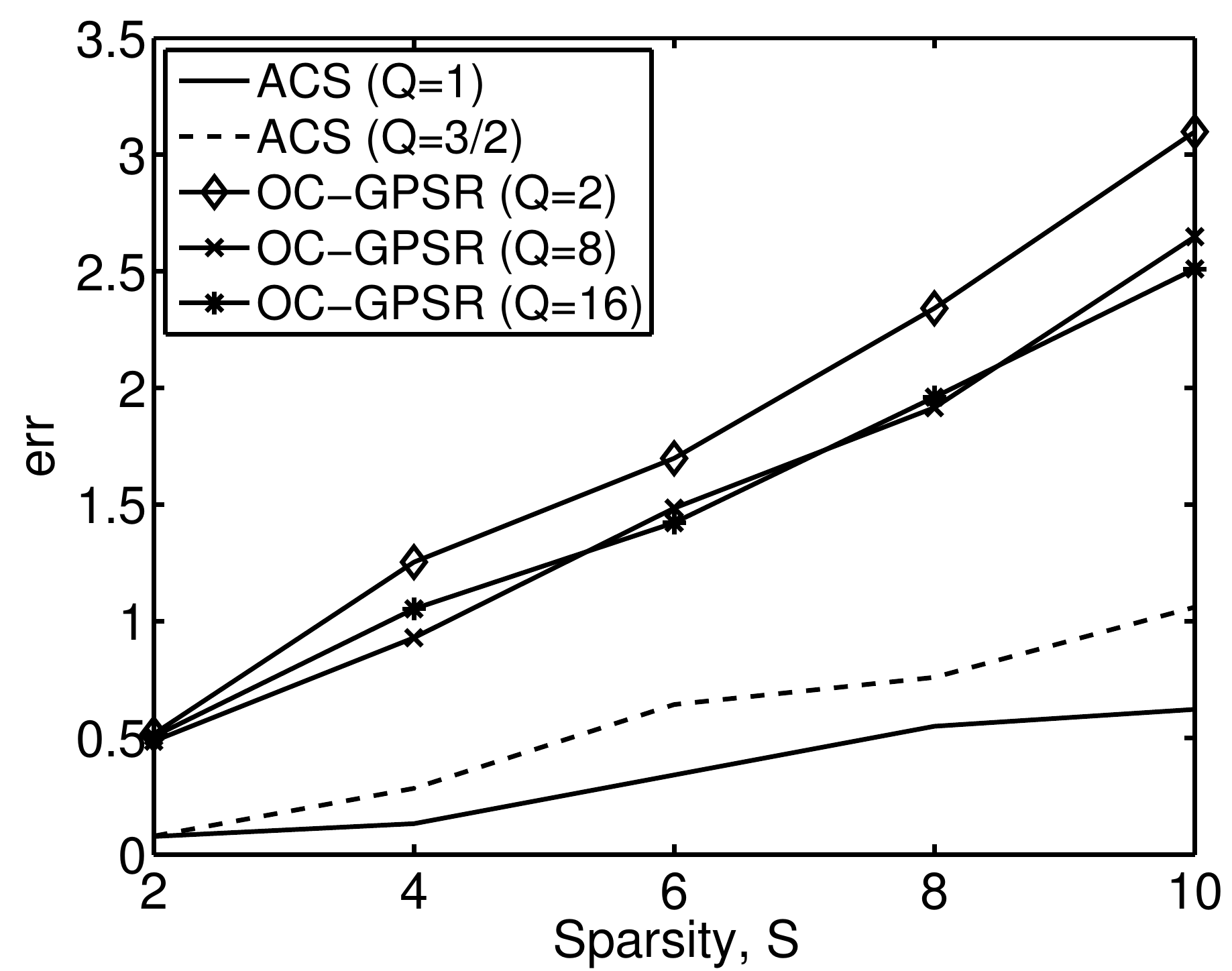}} \\
\subfloat{\includegraphics[width=\figwidth]{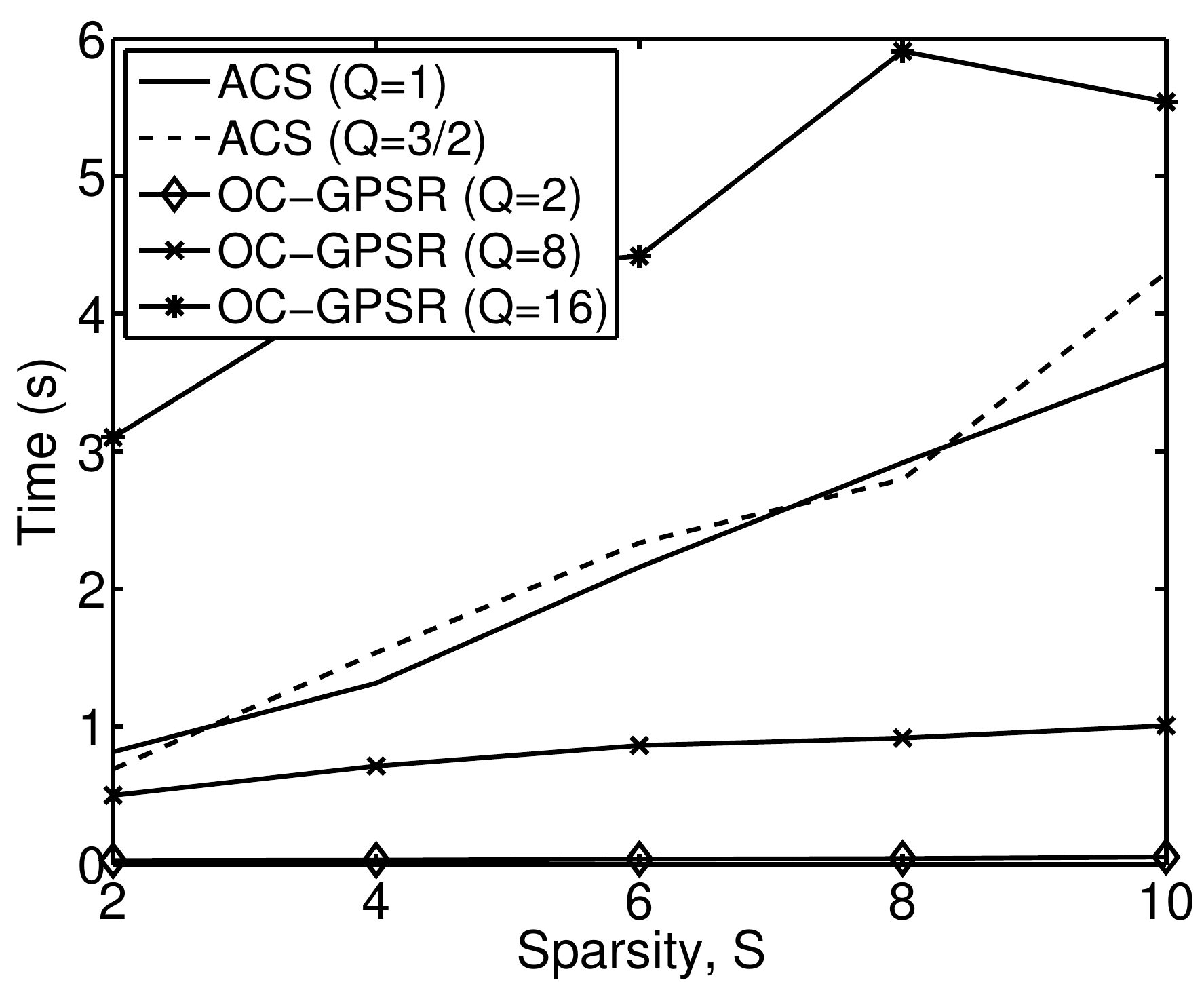}} ~
   \caption{Comparison between the ACS approach and overcomplete dictionaries in their respective ability to recover harmonic signals of varying sparsity.  Shown are the normalized RMSE, the $\mathrm{err}$ metric, and the algorithm runtime. (Averaged over 50 random realizations).}
   \label{fig:oc_comp}
\end{center}
\vspace{-3ex}
\end{figure}

Finally, while the ACS approach is only guaranteed to be biconvex for $Q \ge 3/2$, we see similar RMSE values as for the $Q=1$ case, suggesting that $Q \ge 3/2$ is not necessary in practice and good results can be obtained by perturbing the traditional Fourier basis.   We have found this to be true for a variety of sensing matrices and problem sizes.  In the results of the next section, we therefore use only ACS with $Q=1$.

\subsection{Random Temporal Sampling Strategy}

As a final experiment, we again consider the signal model (\ref{eqn:sigmodel}), however, we use the random temporal sampling strategy used for SDP \cite{Tang:12} and GFB \cite{Rao:12} to form 
the compressed samples.  That is to say, the compressed samples are a random subset of all possible $z(\tau_i)$ elements, without replacement. In this case the 
sensing matrix ${\bf A}$ is a binary matrix consisting of a single ``$1$'' in each row, corresponding to the selected column $\tau_j$.  This is precisely the problem studied in 
\cite{Tang:12} and \cite{Rao:12}, and hence forms a good basis for a comparison of the methods. 
We note that it is for this random sampling architecture that a comparison to the SDP and GFB approaches can be easily made, without having to make nontrivial modifications to the provided codes. On the other hand, ACS easily applies to generic ${\bf A}$.

\setlength{\figwidth}{0.4\linewidth}
\begin{figure}[t]
\begin{center}
\subfloat{\includegraphics[width=\figwidth]{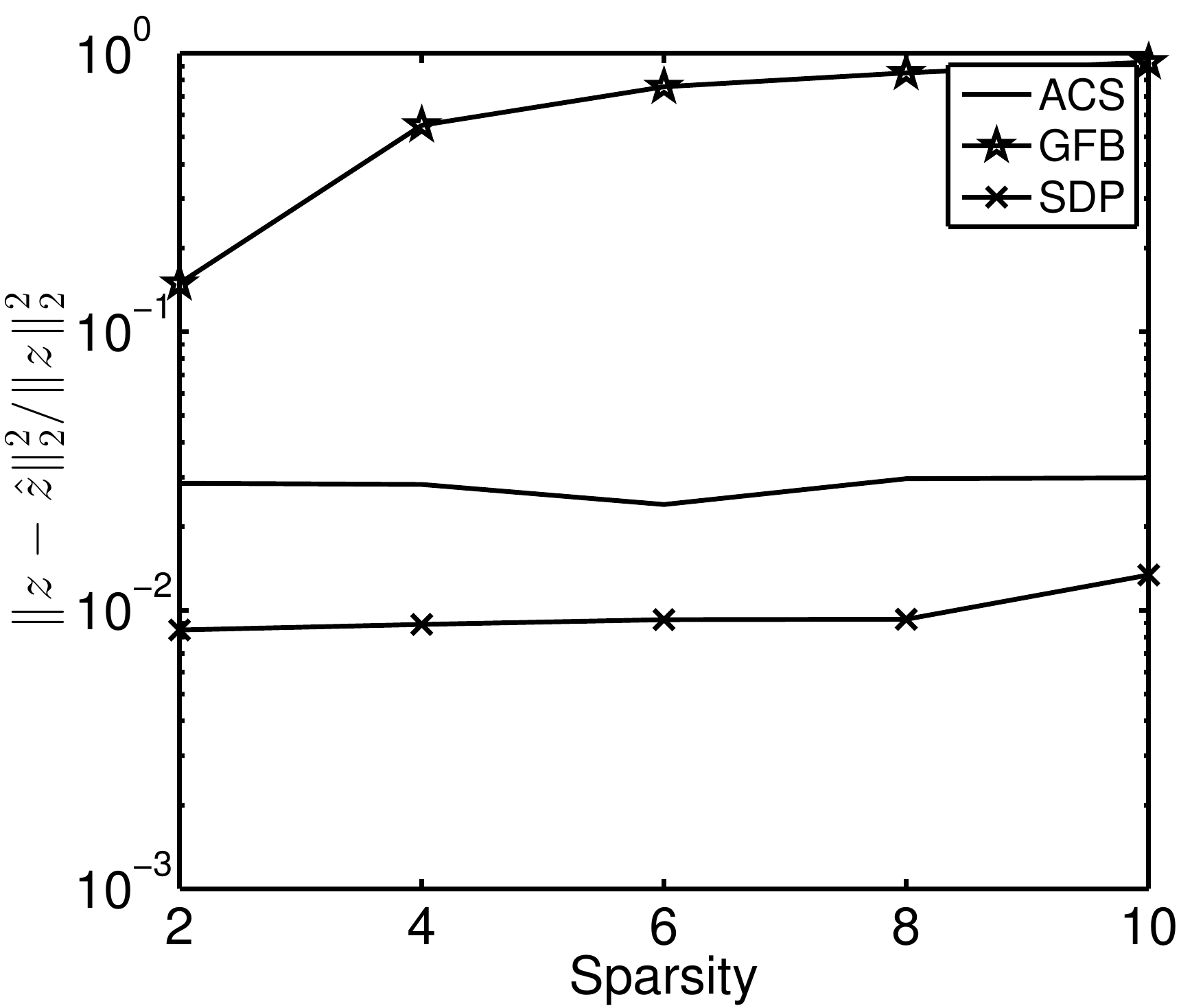}} ~
\subfloat{\includegraphics[width=\figwidth]{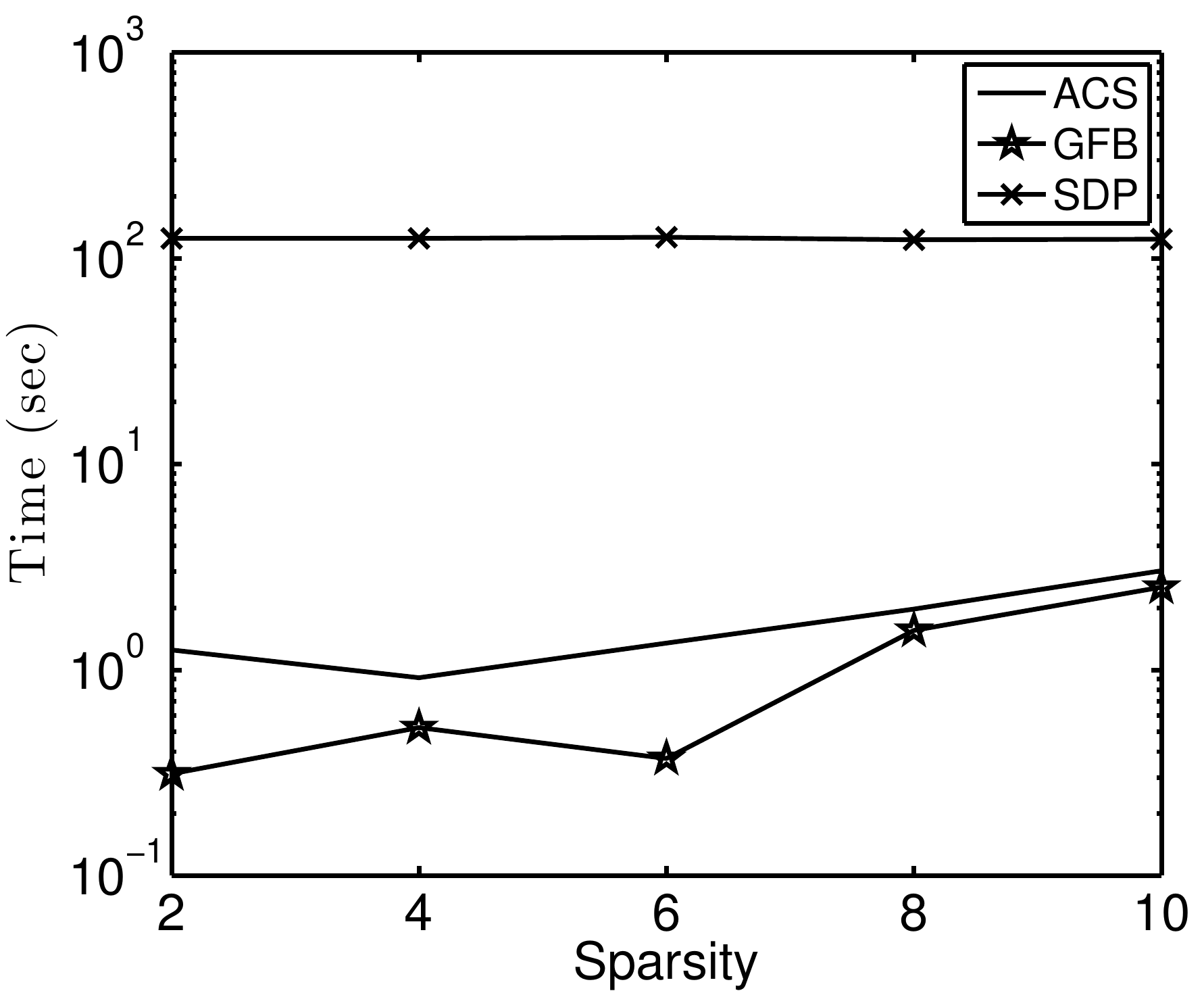}} ~
   \caption{Expected RMSE (left) and runtimes (right) associated with the ACS, SDP, and GFB as a function of signal sparsity. (Averaged over 50 random realizations).}
   \label{fig:avg2}
\end{center}
\vspace{-3ex}
\end{figure}

Comparisons can be seen in Figure \ref{fig:avg2}.
With regard to recovery error, the SDP produces the most accurate recovery with ACS offering comparable performance.  By contrast, the GFB approach tends to produce considerably larger RMSE. However, with regard to computation time, both GFB and ACS are very fast. By contrast, SDP is quite slow, taking roughly 2 minutes for all sparsity levels.

We note that only no-phase signals were considered in this simulation. This was because GFB did not seem robust to phase changes, which was not trivial to fix. Both ACS and SDP, on the other hand, were shown to be robust to phase changes and remain the two main competitors for an arbitrary signal where phases may be unknown. 

In short, ACS offers a good balance between computation time and accuracy.  It runs faster than SDP, while nearly matching it in accuracy.  Moreover, the ACS approach does not require known phases, and is flexible in the sense that it is easily implementable with any valid CS matrix, not just one for random temporal sampling.  ACS also produces better overall performance (RMSE, $\mathrm{err}$) than an overcomplete dictionary.

\section{Conclusions \label{sec:conclusions}}  

This work describes an Alternating Convex Search (ACS) algorithm for correcting frequency errors in the signal model used in compressive sampling 
applications for real harmonic signals.  The algorithm treats both frequencies and model coefficients as unknowns and uses an iterative approach in developing estimates.  
Specifically, the approach uses the familiar GPSR algorithm to update the model coefficients, followed by a maximum likelihood estimate of the unknown 
frequency locations. The algorithm was demonstrated effective at recovering harmonic signals possessing varying levels of sparsity in competitive computation times.
Finally, we note that the code associated with this paper is available to download at: \url{https://sites.google.com/site/albertktoh/software/acs}.

\bibliographystyle{IEEEtran}
\bibliography{refs}

\begin{thebibliography}{10}
\providecommand{\url}[1]{#1}
\csname url@samestyle\endcsname
\providecommand{\newblock}{\relax}
\providecommand{\bibinfo}[2]{#2}
\providecommand{\BIBentrySTDinterwordspacing}{\spaceskip=0pt\relax}
\providecommand{\BIBentryALTinterwordstretchfactor}{4}
\providecommand{\BIBentryALTinterwordspacing}{\spaceskip=\fontdimen2\font plus
\BIBentryALTinterwordstretchfactor\fontdimen3\font minus
  \fontdimen4\font\relax}
\providecommand{\BIBforeignlanguage}[2]{{%
\expandafter\ifx\csname l@#1\endcsname\relax
\typeout{** WARNING: IEEEtran.bst: No hyphenation pattern has been}%
\typeout{** loaded for the language `#1'. Using the pattern for}%
\typeout{** the default language instead.}%
\else
\language=\csname l@#1\endcsname
\fi
#2}}
\providecommand{\BIBdecl}{\relax}
\BIBdecl

\bibitem{Donoho:09}
A.~M. Bruckstein, D.~L. Donoho, and M.~Elad, ``From sparse solutions of systems
  of equations to sparse modeling of signals and images,'' \emph{SIAM Review},
  vol.~51, no.~1, pp. 34--81, 2009.

\bibitem{Donoho:06}
D.~L. Donoho, ``Compressed sensing,'' \emph{IEEE Transactions Information
  Theory}, vol.~52, no.~4, pp. 1289--1306, 2006.

\bibitem{Candes:06a}
E.~J. Cand\`{e}s, J.~Romberg, and T.~Tao, ``Robust uncertainty principles:
  Exact signal reconstruction from highly incomplete frequency information,''
  \emph{IEEE Transactions Information Theory}, vol.~52, no.~2, pp. 489--509,
  2006.

\bibitem{Donoho:10}
D.~L. Donoho and J.~Tanner, ``Precise undersampling theorems,''
  \emph{Proceedings of the IEEE}, vol.~98, no.~6, 2010.

\bibitem{Candes:05}
E.~J. Cand\`{e}s and T.~Tao, ``Decoding by linear programming,'' \emph{IEEE
  Transactions Information Theory}, vol.~15, no.~12, pp. 4203--4215, 2005.

\bibitem{Candes:12}
E.~J. Cand\`{e}s and C.~F.-Granda, ``Towards a mathematical theory of
  super-resolution,'' 2012, arXiv:1203.5871v3.

\bibitem{Tropp:10}
J.~A. Tropp, J.~N. Laska, M.~F. Duarte, J.~K. Romberg, and R.~G. Baraniuk,
  ``Beyond {N}yquist: Efficient sampling of sparse bandlimited signals,''
  \emph{IEEE Transactions on Information Theory}, vol.~56, no.~1, pp. 520--544,
  2010.

\bibitem{Nichols:11}
J.~M. Nichols and F.~Bucholtz, ``Beating {N}yquist with light: A compressively
  sampled photonic link,'' \emph{Optics Express}, vol.~19, no.~8, pp.
  7339--7348, 2011.

\bibitem{Yenduri:12}
P.~K. Yenduri, A.~Z. Rocca, A.~S. Rao, S.~Naraghi, M.~P. Flynn, and A.~C.
  Gilbert, ``A low-power compressive sampling time-based analog-to-digital
  converter,'' \emph{IEEE Journal on Emerging and Selected Topics in Circuits
  and Systems}, vol.~2, no.~3, pp. 502--515, 2012.

\bibitem{Romberg:08}
J.~Romberg, ``Imaging via compressive sampling,'' \emph{IEEE Signal Processing
  Magazine}, vol.~25, no.~2, pp. 14--20, 2008.

\bibitem{lustig2007sparse}
M.~Lustig, D.~Donoho, and J.~M. Pauly, ``Sparse {MRI}: The application of
  compressed sensing for rapid mr imaging,'' \emph{Magnetic resonance in
  medicine}, vol.~58, no.~6, pp. 1182--1195, 2007.

\bibitem{willett2014sparsity}
R.~Willett, M.~Duarte, M.~Davenport, and R.~Baraniuk, ``Sparsity and structure
  in hyperspectral imaging: Sensing, reconstruction, and target detection,''
  \emph{Signal Processing Magazine, IEEE}, vol.~31, no.~1, pp. 116--126, 2014.

\bibitem{harmany2011spatio}
Z.~T. Harmany, R.~F. Marcia, and R.~M. Willett, ``Spatio-temporal compressed
  sensing with coded apertures and keyed exposures,'' 2011, arXiv preprint
  arXiv:1111.7247.

\bibitem{Nichols:13}
C.~V. McLaughlin, J.~M. Nichols, and F.~Bucholtz, ``Basis mismatch in a
  compressively sampled photonic link,'' \emph{IEEE Photonics Technology
  Letters}, vol.~25, no.~23, pp. 2297--2300, December 2013.

\bibitem{Calderbank:11}
Y.~Chi, L.~Scharf, A.~Pezeshki, and A.~R. Calderbank, ``Sensitivity to basis
  mismatch in compressed sensing,'' \emph{IEEE Transactions on Signal
  Processing}, vol.~59, no.~5, pp. 2182--2195, 2011.

\bibitem{Gilbert:12}
P.~Boufounos, V.~Cevher, A.~C. Gilbert, Y.~Li, and M.~J. Strauss, ``What's the
  frequency, {K}enneth?: Sublinear fourier sampling off the grid,''
  \emph{Lecture notes in Computer Science}, vol. 7408, pp. 61--72, 2012.

\bibitem{Tang:12}
G.~Tang, B.~N. Bhaskar, P.~Shah, and B.~Recht, ``Compressive sensing off the
  grid,'' in \emph{Communication, Control, and Computing (Allerton), 2012 50th
  Annual Allerton Conference on}.\hskip 1em plus 0.5em minus 0.4em\relax IEEE,
  2012, pp. 778--785.

\bibitem{Ekanadham:11}
C.~Ekanadham, D.~Tranchina, and E.~P. Simoncelli, ``Recovery of sparse
  translation-invariant signals with continuous basis pursuit,'' \emph{IEEE
  Transactions on Signal Processing}, vol.~59, no.~10, 2011.

\bibitem{Rao:12}
N.~Rao, P.~Shah, S.~Wright, and R.~Nowak, ``A greedy forward-backward algorithm
  for atomic norm constrained minimization,'' in \emph{Acoustics, Speech and
  Signal Processing (ICASSP), 2013 IEEE International Conference on}.\hskip 1em
  plus 0.5em minus 0.4em\relax IEEE, 2013, pp. 5885--5889.

\bibitem{Gorski:07}
J.~Gorski, F.~Pfeuffer, and K.~Klamroth, ``Biconvex sets and optimization with
  biconvex functions: A survey and extensions,'' \emph{Mathematical Methods of
  Operations Research}, vol.~66, pp. 373--407, 2007.

\bibitem{Valley:13}
G.~C. Valley, G.~A. Sefler, and T.~J. Shaw, ``Sensing rf signals with the
  optical wideband converter,'' in \emph{Proceedings of the SPIE, Broadband
  Access Communication Technologies VII}, vol. 8645, San Francisco, CA, USA,
  February 2013, p. doi:10.1117/12.2002144.

\bibitem{Donoho:01}
S.~S. Chen, D.~L. Donoho, and M.~A. Saunders, ``Atomic decomposition by basis
  pursuit,'' \emph{SIAM Review}, vol.~43, no.~1, pp. 129--159, 2001.

\bibitem{Stoica:89}
P.~Stoica, R.~L. Moses, B.~Friedlander, and T.~Soderstrom, ``Maximum likelihood
  estimation of the parameters of multiple sinusoids from noisy measurements,''
  \emph{IEEE Transactions on Acoustics, Speech, and Signal Processing},
  vol.~37, no.~3, pp. 378--392, 1989.

\bibitem{gpsr}
M.~A.~T. Figueiredo, R.~D. Nowak, and S.~J. Wright, ``Gradient projection for
  sparse reconstruction: Application to compressed sensing and other inverse
  problems,'' \emph{IEEE Journal of Selected Topics in Signal Processing:
  Special Issue on Convex Optimization Methods for Signal Processing}, vol.~1,
  no.~4, pp. 586--597, 2007.

\bibitem{Recht:14}
G.~Tang, B.~N. Bhaskar, and B.~Recht, ``Near minimax line spectral
  estimation,'' in \emph{Information Sciences and Systems (CISS), 2013 47th
  Annual Conference on}.\hskip 1em plus 0.5em minus 0.4em\relax IEEE, 2013, pp.
  1--6.

\end{thebibliography}

\clearpage

\begin{appendix} \label{app:A}

\renewcommand{\theequation}{A\arabic{equation}}
\setcounter{equation}{0}  
\setcounter{figure}{0}  
\section*{Convexity of Surrogate Cost Function}  

In this section, we discuss the convexity of the function (\ref{eqn:cost}).  This is the primary requirement for making statements about algorithm convergence.  In particular, we may make use of the following theorem:
\begin{theorem}
(Theorem 4.5 in \cite{Gorski:07}) Let $B\subset\R^{M}\times\R^{N}$, let $f:B\rightarrow\R$ be a biconvex function which is bounded from below, and let the optimization problems (\ref{eqn:cost1},\ref{eqn:cost2}) be 
solvable.  Then the sequence $f(\{ {\bf x}^{(i)},{\bm\theta^{(i)}} \})_{i\in \N}$ generated by ACS converges monotonically \cite{Gorski:07}.
\label{theorem:converge}
\end{theorem}
For fixed ${\bm\theta}^*$ the problem (\ref{eqn:cost1}) is known to be convex in ${\bf x}$. 

While the function described by (\ref{eqn:cost2}) is not convex, it can be well-approximated by a function that is convex provide each of the frequency parameter is restricted to the range $\theta_k\in[-1/(2QN),1/(2QN)]$, and where $Q$ is determined in what follows.  To begin, denote the true and estimated signals at a given iteration of the algorithm as ${\bf z}$, where each element is $z_n=\sum_k x_k\sin(2\pi f_k n), n=0\cdots N-1,~f_k\in[0,1/2]$ and $\hat{\bf z}({\bm\theta})$, where each element is $\hat{z}({\bm\theta})_n=\sum_k \hat{x}_k\sin(2\pi (f_k+\theta_k)n), n=0\cdots N-1$ respectively.  We are interested in the cost function
\begin{align}
\|{\bf y}-{\bf A}\hat{\bf z}\|_2^2&=\|{\bf Az}-{\bf A}\hat{\bf z}({\bm\theta})+{\bm \eta}\|_2^2\nonumber \\
&=\|{\bf A}({\bf z}-\hat{\bf z}({\bm\theta}))+{\bm \eta}\|_2^2\nonumber \\
&\equiv \|{\bf Av}_{\bm\theta}+{\bm \eta}\|_2^2,
\label{eqn:app1}
\end{align}
where ${\bf A}$ is an $M\times N$ matrix comprised of i.i.d. entries chosen from some probability distribution 
function with zero mean and variance $\sigma_A^2$.  The noise vector ${\bm \eta}$ is an $M\times 1$ vector consisting of 
zero-mean, i.i.d. Gaussian entries $\eta_m$, each with variance $\sigma^2$. The $N\times 1$ vector ${\bf v}_{\bm\theta}$ is a measure of the discrepancy between signal and estimate where the subscript signifies this vector is a function of our unknown parameters ${\bm\theta}$.  Each member of this vector will be denoted $v_{n\bm\theta},~n=1\cdots N$.  Similarly, each member of ${\bf A}$ will be denoted $a_{mn}$ and each projection of the difference vector will be written ${\bf a}_m{\bf v}_{\bm\theta}$ where we have denoted ${\bf a}_m$ as the $m^{th}$ row of the matrix ${\bf A}$.  Moreover, it should be noted that each $E[a_{mn}v_{n\bm\theta}]=0$, and $E[(a_{mn}v_{n\bm\theta})^2]=\sigma_A^2 v^2_{n\bm\theta}$ where the expectation is taken with respect to the probability distribution governing the $a_{mn}$.  

Expanding (\ref{eqn:app1}) and noting that the minimum of the cost function is not affected by re-scaling, we have
\begin{align}
\frac{1}{M}\|{\bf Av}_{\bm\theta}+{\bm \eta}\|_2^2&=\frac{1}{M}\left\{{\bf v}_{\bm\theta}^T{\bf A}^T{\bf Av}_{\bm\theta}+2{\bm \eta}^T{\bf Av}_{\bm\theta}+{\bm \eta}^T{\bm \eta}\right\}\nonumber \\
\label{eqn:expand}
\end{align}
In what follows we will make use of the central limit theorem in approximating this cost function.  Recall that the central limit theorem states that for $M$ independent realizations, $x_m$, of a random variable $X$ with mean $\mu_X$ and variance $\sigma_X^2$, the summation $\frac{1}{M}\sum_{m}x_m\sim \mathcal{N}(\mu_X,\sigma_X^2/M)$ for $M$ large.  In other words, the discrete summation is equal to the true value plus an error term of order $O(\sigma_X/\sqrt{M})$.  

With this in mind, we first note that the last term is simply the noise variance
\begin{align}
\frac{1}{M}\sum_{m=0}^{M-1} \eta_m^2=\sigma^2+O(\sqrt{2}\sigma^2/\sqrt{M})
\end{align}
and provides a constant offset to the cost function.  The order of the error is determined by noting that for a normally distributed random variable $X$, with variance $\sigma_X^2$, the variance of $X^2$ is $2\sigma_X^4$.  

The cost function is largely governed by ${\bf v}_{\bm\theta}^T{\bf A}^T{\bf Av}_{\bm\theta}$ which is a non-convex function of the unknown parameter vector ${\bm\theta}\in\R^N$. However, over a convex set $\Theta$ which restricts the range of each $\theta_k$, this function can be well approximated by a function that is convex.   To show this, we first re-write this term as
\begin{align}
\frac{1}{M}{\bf v}_{\bm\theta}^T&{\bf A}^T{\bf Av}_{\bm\theta}=\frac{1}{M}\sum_{m=0}^{M-1}({\bf a}_m{\bf v}_{\bm\theta})^2=\frac{1}{M}\sum_{m=0}^{M-1}\left(\sum_{n=0}^{N-1} a_{mn}v_{n\bm\theta}\right)^2\nonumber \\
&=\frac{1}{M}\sum_{m=0}^{M-1} \left[\sum_{n=0}^{N-1} (a_{mn}v_{n\bm\theta})^2+\sum_{n=0}^{N-1} \sum_{q\ne n}^Na_{mn}a_{mq} v_{n\bm\theta}v_{q\bm\theta}\right]\nonumber \\
&=\sum_{n=0}^{N-1} v_{n\bm\theta}^2 \frac{1}{M}\sum_{m=0}^{M-1} a_{mn}^2+\sum_{n=0}^{N-1} \sum_{q\ne n}^Nv_{n\bm\theta}v_{q\bm\theta}\frac{1}{M}\sum_{m=0}^{M-1} a_{mn}a_{mq} \nonumber \\
&=\sum_{n=0}^{N-1} v_{n\bm\theta}^2 \left\{E\left[a_{mn}^2\right]+O\left(\frac{\sqrt{2+\gamma(A)}\sigma_A^2}{\sqrt{M}}\right)\right\}\nonumber \\
&~~+\sum_{n=0}^{N-1} \sum_{q\ne n}^N v_{n\bm\theta}v_{q\bm\theta} \left\{E\left[a_{mn}a_{qn}\right]+O\left(\frac{\sigma_A^2}{\sqrt{M}}\right)\right\}\nonumber \\
&=\sigma_A^2 \sum_{n=0}^{N-1} v_{n\bm\theta}^2\nonumber \\
&~~ + O\left(\frac{\sqrt{2+\gamma(A)}\sigma_A^2}{\sqrt{M}}\sum_{n=0}^{N-1} v_{n\bm\theta}^2+\frac{\sigma_A^2}{\sqrt{M}}\sum_{n=0}^{N-1} \sum_{q\ne n}^N v_{n\bm\theta}v_{q\bm\theta} \right).
\label{eqn:cross}
\end{align}
The constant $\gamma(X)\equiv \frac{E[(X-\mu)^4]}{\sigma_X^4}-3$ is the kurtosis associated with the random variable $X$ and depends on the probability distribution governing $X$.  For example, if the $a_{mn}$ are normally distributed $\gamma(A)=0$ while for entries chosen from a uniform distribution, $\gamma(A)=-6/5$.   Furthermore we have under the assumption that the values $a_{mn}$ are independently chosen, $E[a_{mn}a_{qn}]=0$ while the variance of the second term in (\ref{eqn:cross}) is the product of the variances of the constituent variables.  Thus, by the central limit theorem, approximating this covariance term with an expectation yields an error of order $\sigma_A^4$.

In short, the degree to which the cost function is approximated by the last line in (\ref{eqn:cross}) is entirely predicated on the behavior of the function 
\begin{align}
\sum_{n=0}^{N-1}&v_{n\bm\theta}^2\nonumber \\
&=\sum_{n=0}^{N-1}\left[\sum_{k=1}^S x_k\sin(2\pi f_k n)-\sum_{k=1}^S 
\hat{x}_k\sin(2\pi (f_k+\theta_k)n)\right]^2\nonumber \\
&=\sum_{n=0}^{N-1} \left[\sum_{k=1}^S \left[x_k\sin(2\pi f_kn)-\hat{x}_k\sin(2\pi(f_k+\theta_k)n)\right]^2\right.\nonumber \\
&\left.+\sum_{k=1}^S\sum_{k\ne k'}x_kx_{k'}
\sin(2\pi f_kn)\sin(2\pi f_{k'}n)\right.\nonumber \\
&\left.-2\sum_{k=1}^S\sum_{k\ne k'} x_k\hat{x}_k\sin(2\pi f_k n)\sin(2\pi f_{k'}n)\right.\nonumber \\
&\left.+\sum_{k=1}^S\sum_{k\ne k'}\hat{x}_k\hat{x}_{k'}\sin(2\pi 
f_kn)\sin(2\pi (f_k+\theta_k)n)\right].
\end{align}
The asymptotic (large $N$) properties of the above are easily obtained by noting that
\begin{align}
\lim_{N\rightarrow\infty} \sum_{n=0}^{N-1}\sin(2\pi f_k n)\sin(2\pi f_{k'}n)= 0.
\end{align}
The result is that
\begin{align}
\lim_{N\rightarrow\infty} &\sum_{n=0}^{N-1}v_{n\bm\theta}^2\nonumber \\
&=\lim_{N\rightarrow\infty}\sum_{k=1}^{N}\sum_{n=0}^{N-1} \left[x_k\sin(2\pi f_k n)-\hat{x}_k\sin(2\pi (f_k+\theta_k) n)\right]^2\nonumber \\
&=\lim_{N\rightarrow\infty} \sum_{k=1}^{N} g_k(\theta_k)
\end{align}
where 
\begin{align}
g_k(\theta_k)&=\sum_{n=0}^{N-1} \left[x_k\sin(2\pi f_k n)-\hat{x}_k\sin(2\pi (f_k+\theta_k) n)\right]^2
\label{eqn:app3}
\end{align}
Thus, the $N$-dimensional cost can be approximated by a sum of $N$ one-dimensional functions.  We also have that $g_k(\theta_k)\ge 0$, ensuring that the $N$-dimensional minimum of this function can be found by taking each solution component (each $\theta_k$) as the minimizers of (\ref{eqn:app3}).  

Returning to Eqn. (\ref{eqn:cross}) it can be seen that for large $N$ the second error term vanishes due to the near orthogonality of the product of sinusoids at different frequencies.  Hence, for even moderately large $N$ we have that
\begin{align}
{\bf v}_{\bm\theta}^T{\bf A}^T{\bf Av}_{\bm\theta}&=\sigma_A^2 \sum_{k=1}^S g_k(\theta_k) + O\left(\frac{\sqrt{2+\gamma(A)}\sigma_A^2}{\sqrt{M}}\sum_{k=1}^S g_k(\theta_k)\right).
\end{align}

Finally, the second term in (\ref{eqn:expand}) can be written
\begin{align}
\frac{2}{M}\sum_{m=0}^{M-1} ({\bf a}_m{\bf v}_{\bm\theta})\eta_m&=\frac{2}{M}\sum_{m=0}^{M-1}\sum_{n=0}^{N-1} a_{mn}\eta_mv_{n{\bm\theta}}\nonumber \\
&=2\sum_{n=0}^{N-1} v_{n{\bm\theta}} \left\{E[a_{mn}]E[\eta_m]+O\left(\frac{\sigma_A\sigma}{\sqrt{M}}\right)\right\}\nonumber \\
&=2\sum_{n=0}^{N-1} v_{n{\bm\theta}} \times O\left(\frac{\sigma_A\sigma}{\sqrt{M}}\right)\nonumber \\
&=2\sum_{k=1}^S h(\theta_k) \times O\left(\frac{\sigma_A\sigma}{\sqrt{M}}\right)
\end{align}
where we have noted that for independent random variables the expectation of the product factors.  We have also made the same orthogonally argument as above in order to write the cost as a function of 
\begin{align}
h_k(\theta_k)=\sum_{n=0}^{N-1} \left[x_k\sin(2\pi f_k n)-\hat{x}_k\sin(2\pi (f_k+\theta_k) n)\right].
\end{align}
Thus, for either zero mean noise, or zero mean entries in the matrix ${\bf A}$, this part of the cost function reduces to an error term of order $\frac{2\sigma_A\sigma}{\sqrt{M}}\sum_{k=1}^S h_k(\theta_k)$.  

Incorporating all components of the cost function we therefore have
\begin{align}
\frac{1}{M}&\|{y-\bf A}\hat{\bf z}_{\bm\theta}\|_2^2= 
\sigma_A^2\sum_{k=1}^{S} g_k(\theta_k)+\sigma^2\nonumber \\
&+O\left(\frac{\sqrt{2+\gamma(A)}\sigma_A^2}{\sqrt{M}}\sum_{k=1}^{S} g_k(\theta_k)+\frac{\sigma_A\sigma}{\sqrt{M}}\sum_{k=1}^S h_k(\theta_k)+\frac{\sigma^2}{\sqrt{M}}\right)
\label{eqn:app2}
\end{align}
so that the quality of the approximation is seen to be dependent on the functions $g_k(\theta),~h_k(\theta)$.    By inspection the function $h_k(\theta)$ is simply a sinusoid with offset governed by $x_k$ and an amplitude governed by $\hat{x}_k$ hence, the order of this error term is at most $O(|x_k+\hat{x}_k|)$ for any value of $\theta_k$ and for any number of data $N$.  As will be shown next, the function $h_k(\theta)$ is $O(1/N)$ relative to the term $g_k(\theta)$ and can safely be neglected in the analysis.  

The approximate cost function and the quality of the approximation are therefore determined by the behavior of $g_k(\theta_k)$.  In describing the behavior of this function we also define the difference between the true amplitude and estimated amplitude as $\delta=|x_k-\hat{x}_k|$.  In this case it can be shown that the function possesses a minimum value of $g_k(0)=\frac{N}{2}\delta^2$ and grows to approach a value of $N|\hat{x}_kx_k|+\frac{N}{2}\delta^2$ for $|\theta_k|\ge 1/2/N$.
\begin{figure}[tbc]
  \centerline{
   \begin{tabular}{c}   
        \includegraphics[width=3.0in,angle=0]{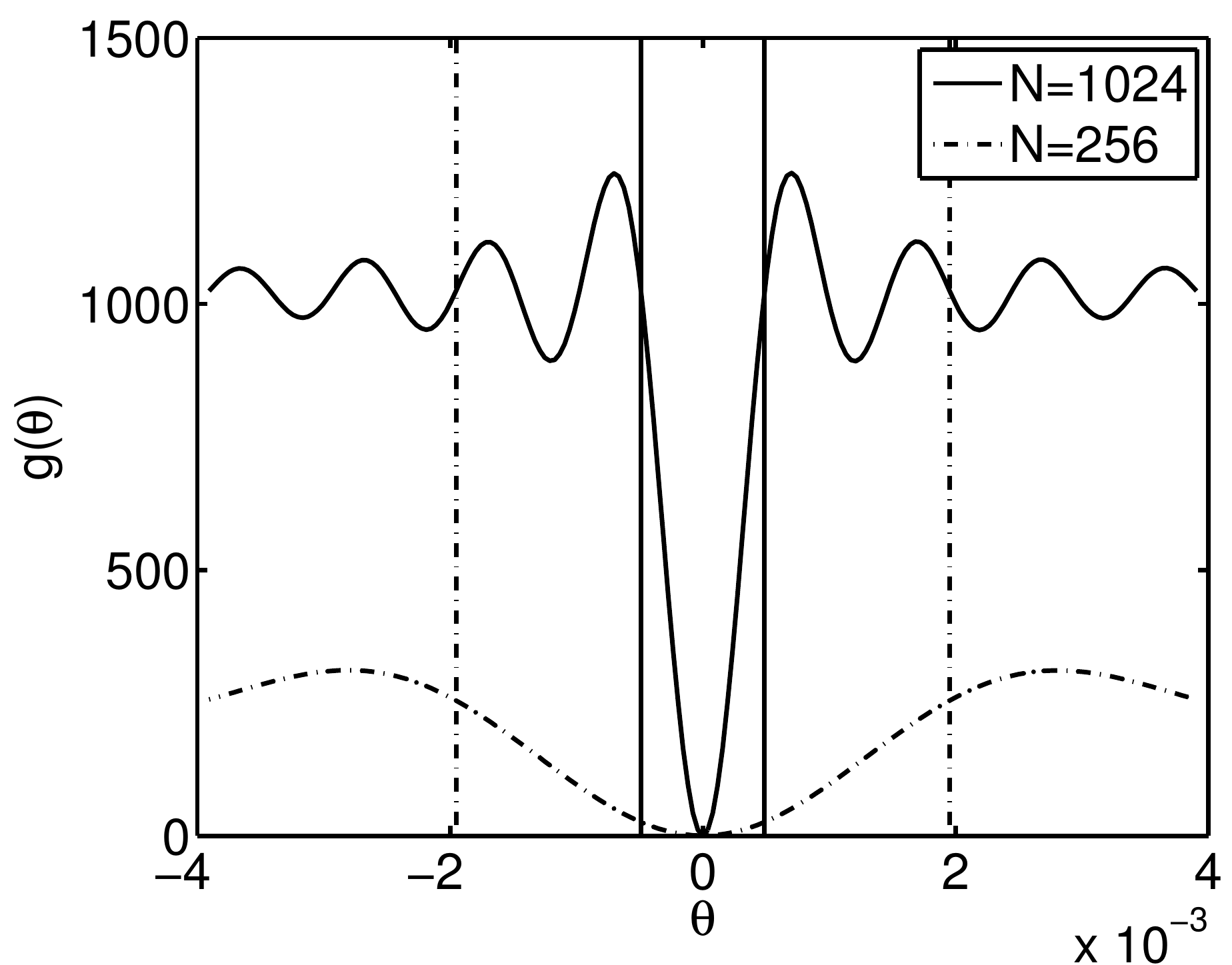}
   \end{tabular}
   }
   \caption{Behavior of the function $g_k(\theta)$.  Also shown are the limits of a standard frequency bin width of $\pm 1/2/N$.  We are ultimately interested in the behavior of this function on this interval for a given value of $N$.}
   \label{fig:H}
\end{figure}
This function is shown in Figure (\ref{fig:H}) for a signal consisting of a single frequency and for two different sets of signal amplitudes.  In the first example we set the true amplitudes equal to the values $\hat{x}=x=1$, hence $\delta=0$.  In the second we set $\hat{x}_k=4,~x_k=1$ to simulate a very large error in the estimated coefficient value.  The approximation is clearly better when the estimated and true amplitudes are close.  However, even with a large estimated coefficient error the function $g_k(0)$ is still the minimum and the shape of the function is very much the same as the true cost function.  In fact, the approximation remains good for $\theta_k$ extending well beyond a standard Fourier bin width.   

Neglecting the higher-order error terms from (\ref{eqn:app2}), we have approximately that
\begin{align}
\frac{1}{M}&\|{y-\bf A}\hat{\bf z}_{\bm\theta}\|_2^2\approx \sigma^2+\sigma_A^2\sum_{k=1}^{S} g_k(\theta_k)+\nonumber \\
&O\left(\frac{\sqrt{2+\gamma(A)}\sigma_A^2SN(|\hat{x}{x}|+\frac{\delta^2}{2})}{\sqrt{M}}\right), \left[-\frac{1}{2N}\le\theta_k\le \frac{1}{2N}\right]\nonumber \\
&\approx \sigma_A^2\sum_{k=1}^{S} g_k(\theta_k)
\label{eqn:appnew}
\end{align}
where $S$ is the number of tones in the signal.  In short, to understand the properties of the stochastic, non-convex function (\ref{eqn:app1}), we may study the properties of the deterministic function (\ref{eqn:appnew}).  

\begin{figure}[tbc]
  \centerline{
   \begin{tabular}{ccc}   
        \includegraphics[width=1.75in,angle=0]{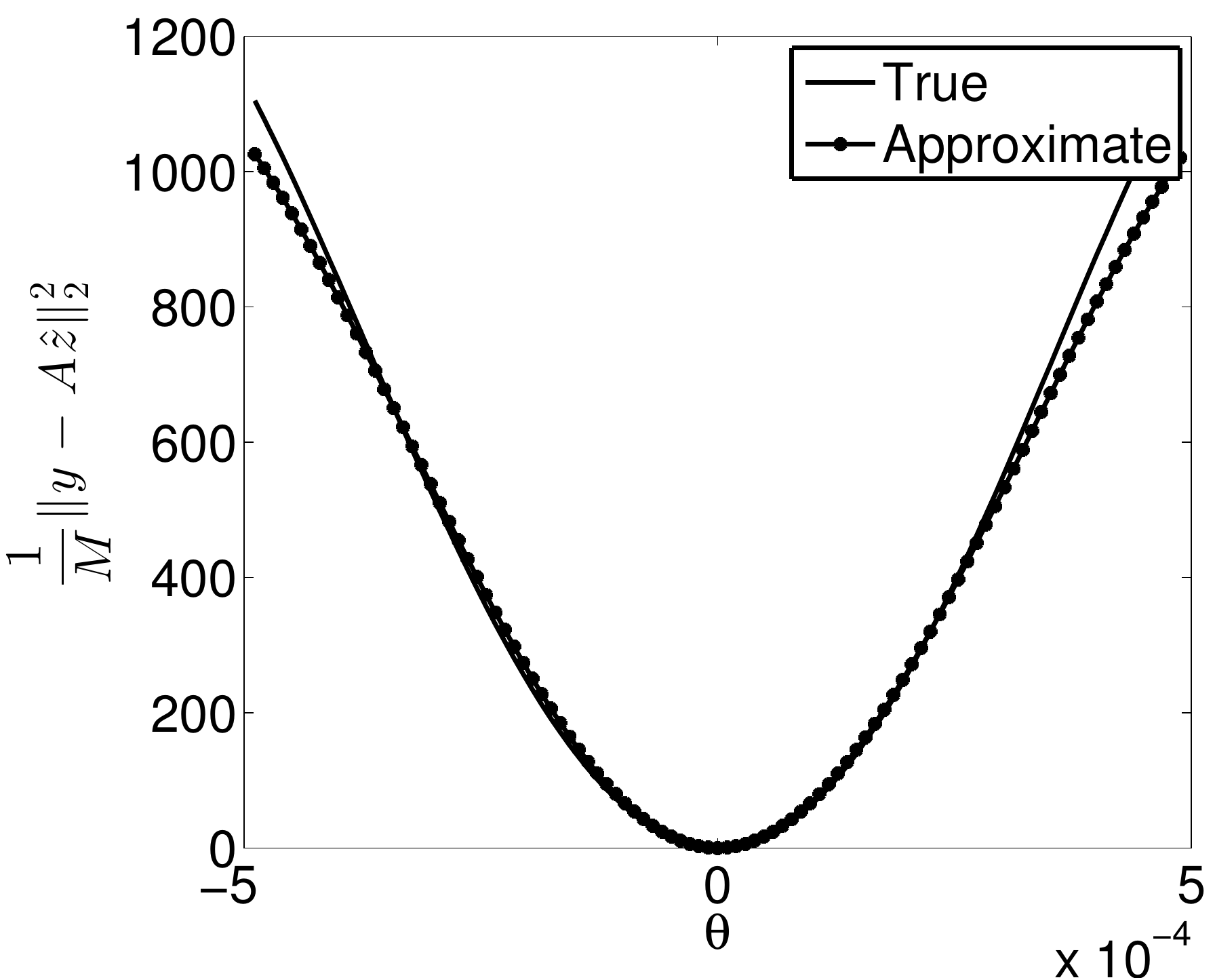} &\hspace*{0.1in} &
         \includegraphics[width=1.75in,angle=0]{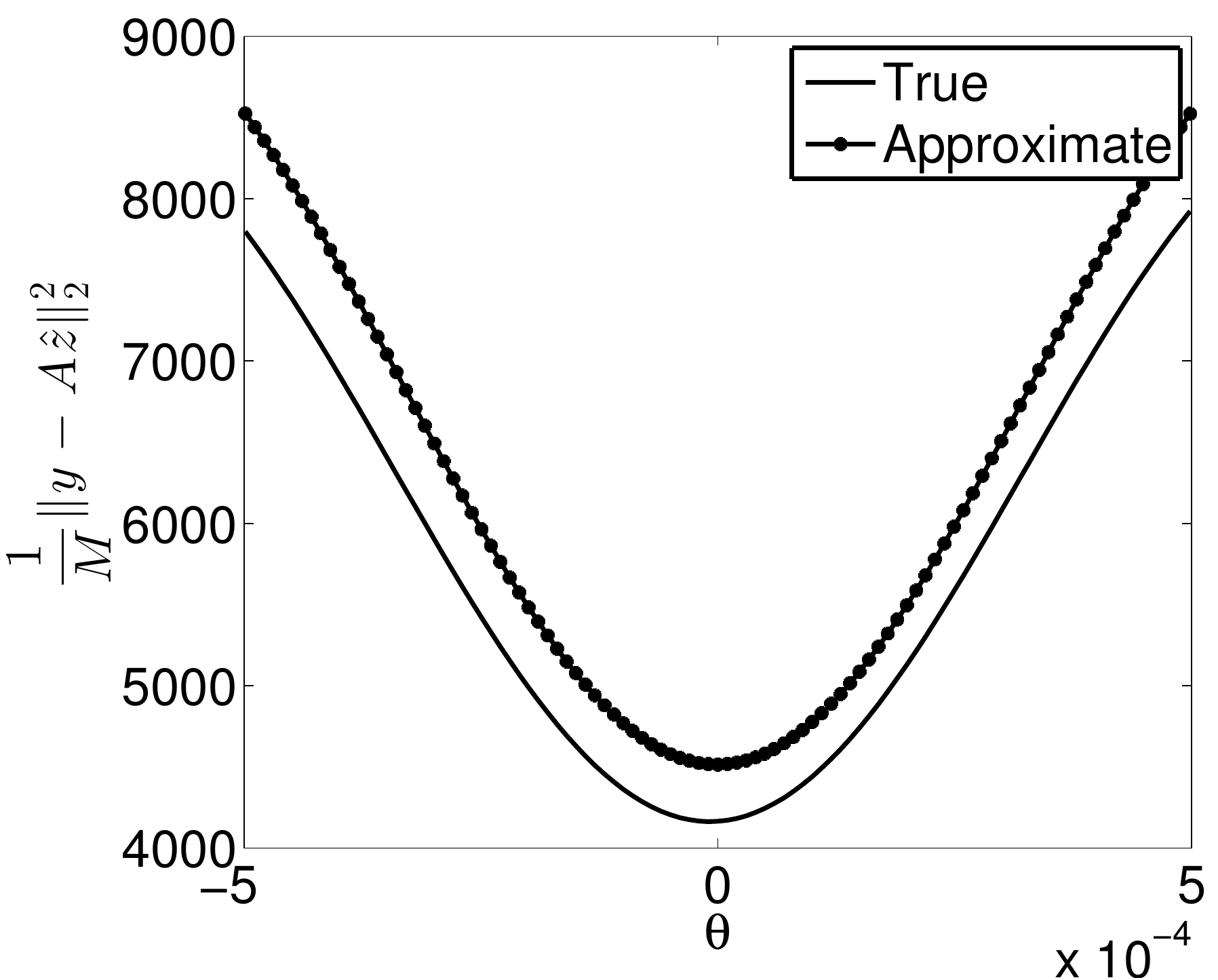}
   \end{tabular}
   }
   \caption{Comparison of cost functions (\ref{eqn:app1}) and (\ref{eqn:app2}) for $N=1024,~M=512$ over a standard frequency bin width, $\theta=[-1/2/N,1/2/N]$.  The entries of $A$ were chosen from a standard normal distribution ($\sigma_A^2=1$).  In the first example (left) we set $\hat{x}=x=1$ so that $\delta=0$.  As predicted, the approximation is best in the vicinity of the true frequency ($\theta=0$), rising to be of order $\sqrt{2}N/\sqrt{M}=63$ near the edges of the bin.  For the second example (right) we set $\hat{x}_1=1$, and $x_1=4$, hence the minimum of the cost function attains the value $\frac{N}{2}\delta^2=4608$ while the error is predicted to be of order $\sqrt{2}\times 8.5N/\sqrt{M}=500$.  Although the discrepancy in estimated and true amplitude drives the cost function to larger values, the minimum remains at the true frequency location, i.e., where $\theta=0$.}
   \label{fig:costcompare}
\end{figure}

As an illustration of the approximate cost (\ref{eqn:appnew}), we choose each $a_{mn}\sim\mathcal{N}(0,1)$, take $M=N=1024$, and plot both (\ref{eqn:app1}) and (\ref{eqn:appnew}) for a single tone ($S=1$) with true frequency $f=0.2$.  We further assumed the case of zero additive noise, i.e., $\sigma=0$.   This comparison is shown in Figure \ref{fig:costcompare}.
We have found the approximation to be accurate to the stated order for a variety of probability distributions governing the entries of ${\bf A}$ and for a variety of sampling parameters $N,~M$.

Given this approximation, we are interested in the range over which (\ref{eqn:appnew}) is convex.  In order to show this we need only demonstrate that each of the $g_k(\theta_k)$ is convex.
Differentiating (\ref{eqn:app3}) twice w.r.t. $\theta_k$, we seek the range over which the function $g''(\theta_k)$ is positive.  Under the large sample approximation we see that the zeros of $g''(\theta_k)$ are given by the zeros of
\begin{align}
g_Z''(\theta_k)&=\lim_{N\rightarrow\infty} g''(\theta_k)=2 N^2 \csc (\pi  \theta_k) \sin (\pi  (2 N+1) \theta_k)\nonumber \\
&+2 N \csc ^2(\pi \theta_k) \cos (2 \pi  (N+1) \theta_k)\nonumber \\
&+3 \csc ^2(\pi \theta_k) \cos (2 \pi  (N+1) \theta_k)\nonumber \\
&+\csc ^3(\pi \theta_k) (-\sin (\pi  (2 N+3) \theta_k))\nonumber \\
&-\cot ^2(\pi  \theta_k)+\csc ^2(\pi  \theta_k)
\label{eqn:secondderiv}
\end{align} 
which is independent of the true and estimated amplitudes $x,~\hat{x}$.  In the discrete setting we are interested in determining the frequency bin width over which this expression is guaranteed to be positive.  In the worst case, each unknown frequency is precisely at the midpoint of each bin in our discretized frequency space.  We therefore re-write the continuous parameter $\theta_k$ in terms of a discrete 
frequency bin size as $\theta_k=\pm \frac{1}{2QN}$ and seek the value of $Q$ for which (\ref{eqn:secondderiv}) is positive.  This substitution yields the transcendental equation
\begin{align}
\tan\left(\frac{\pi}{|Q|}\right)&=\frac{4|Q|\pi}{4|Q|^2-2\pi^2}
\end{align}
for which we require the roots for $|Q|>1$.  It turns out that the only such root is $|Q|=1.509$, so the cost function (\ref{eqn:appnew}) is convex within $\pm \frac{1}{3.018N}$ of the true frequency $f$.  Hence, in order to ensure that the optimization problem being solved is indeed biconvex and that Theorem (\ref{theorem:converge}) applies, we are required to ``tile'' this space with bins of width $\frac{2}{3N}$ which is tantamount to using a Fourier dictionary that is $\frac{3}{2}$ times overcomplete.  In practice, however, we have found that a standard discrete Fourier basis ($Q=1$) is sufficient to produce good convergence.

\end{appendix}

\ifthenelse{\boolean{FIGSINBACK}}{

}

\end{document}